\scrollmode
\documentclass[12pt]{amsart}
\usepackage{verbatim}
\usepackage{eucal}
\usepackage{amssymb}
\usepackage{graphicx}
\usepackage{psfrag}
\usepackage{mathrsfs}
\usepackage{xypic}
\CompileMatrices
\newif \ifwide
\newif \ifavnermargin
\def \makemargins{
\ifwide
	\oddsidemargin .25in
	\evensidemargin .25in
	\textwidth 6.00in
\else
\fi
\ifavnermargin
	\headheight=7pt
	\textheight=574pt
	\textwidth=432pt
	\topmargin=14pt
	\oddsidemargin=18pt
	\evensidemargin=18pt
\else	
\fi
}

\avnermargintrue
\makemargins
\swapnumbers
\theoremstyle{plain}
\newtheorem{theorem}[subsection]{Theorem}
\newtheorem{proposition}[subsection]{Proposition}
\newtheorem{lemma}[subsection]{Lemma}
\newtheorem{corollary}[subsection]{Corollary}

\theoremstyle{definition}
\newtheorem{definition}[subsection]{Definition}
\newtheorem{example}[subsection]{Example}

\theoremstyle{remark}
\newtheorem{remark}[subsection]{Remark}

\newcount\timehh\newcount\timemm
\timehh=\time 
\divide\timehh by 60 \timemm=\time
\newcommand{\draftauthor}[1]{\author{#1
    {
      --- \protect \protect\sc\today\ ---
      \ifnum\timehh<10 0\fi\number\timehh\,:\,\ifnum\timemm<10 0\fi\number\timemm
      \protect \, \, \protect \bf DRAFT
    }
  }
}
\newcommand{\R}{{\mathbb R}}
\newcommand{\A}{{\mathscr A}}

\newcommand{\Z}{{\mathbb Z}}
\newcommand{\Q}{{\mathbb Q}}

\newcommand{\one}{{\mathbf{1}}}

\newcommand{\B}{{\mathscr B}}
\newcommand{\BB}{{\mathbb B}}

\renewcommand{\O}{{\mathscr O}}

\newcommand{\fff}{{\mathfrak f}}
\newcommand{\bbb}{{\mathfrak b}}
\newcommand{\aaa}{{\mathfrak a}}
\newcommand{\SL}{{\mathfrak s \mathfrak l}}
\renewcommand{\ggg}{{\mathfrak g}}
\renewcommand{\hat}{\widehat}

\newcommand{\eee}{{{\mathbf{e}}}}

\newcommand{\sets}[1]{[\![#1]\!]}
\newcommand{\Index}[1]{\|#1\|}

\DeclareMathOperator*{\Max}{Max}

\DeclareMathOperator{\Diag}{diag}
\DeclareMathOperator{\Sign}{sgn}

\DeclareMathOperator{\Rank}{rank}

\DeclareMathOperator{\Id}{Id}
\DeclareMathOperator{\Trace}{Tr}
\begin{document}

\title[Evaluation of Dedekind sums]{Evaluation of Dedekind sums,
Eisenstein cocycles, and special values of $L$-functions}

\newif \ifdraft
\def \makeauthor{
\ifdraft
	\draftauthor{Paul E. Gunnells and Robert Sczech}
\else
\author{Paul E. Gunnells}
\address{Department of Mathematics\\
Columbia University\\
New York, NY  10027}
\email{gunnells@math.columbia.edu}

\author{Robert Sczech}
\address{Department of Mathematics\\
Rutgers University--Newark\\
Newark, NJ  07102}
\email{sczech@andromeda.rutgers.edu}

\fi
}

\thanks{The first author was partially supported by a Columbia
University Faculty Research grant and the NSF}

\draftfalse
\makeauthor

\ifdraft
	\date{\today}
\else
	\date{September 23, 1999}
\fi



%
%

\begin{abstract}
We define certain higher-dimensional Dedekind sums that generalize the
classical Dedekind-Rademacher sums, and show how to compute them
effectively using a generalization of the continued-fraction
algorithm.  

We present two applications.  First, we show how to express special
values of partial zeta functions associated to
totally real number fields in terms of these sums via the Eisenstein
cocycle introduced by the second author.  Hence we obtain a
polynomial-time algorithm for computing these special values.  Second,
we show how to use our techniques to compute certain special values of
the Witten zeta-function, and compute some explicit examples.
\end{abstract}
\maketitle

%
%
\section{Introduction}\label{introduction}

\subsection{}\label{intro}
Let $\sigma $ be a square matrix with integral columns $\sigma
_{j}\in \Z ^{n}$ ($j=1,\dots ,n$), and let $L\subset \Z ^{n}$ be a
lattice of rank $\geq 1$.  Let $v\in \R^{n}$, and let $e\in \Z ^{n}$
with $e_{j}\geq 1$.  Associated to the data $(L,\sigma ,e,v)$ is the
\emph{Dedekind sum}
\begin{equation}\label{Dede.sum}
S = S (L,\sigma ,e,v) := \sideset{}{'}\sum _{x\in L} \eee (\langle
x,v\rangle )\frac{\det \sigma}{\langle x,\sigma _{1}\rangle
^{e_{1}}\cdots \langle x,\sigma _{n}\rangle ^{e_{n}}}.
\end{equation}
Here $\langle x,y\rangle := \sum x_{i}y_{i}$ is the usual scalar
product on $\R ^{n}$, $\eee (t)$ is the character $\exp (2\pi it)$,
and the prime next to the summation means to omit terms for which the
denominator vanishes.  This series converges absolutely if all
$e_{j}>1$, but only conditionally if $e_{j}=1$ for some $j$.  In
the latter case, we define the value of $S$ by the $Q$-limit
\begin{equation}
\sideset{}{'}\sum _{x\in L} a (x)\Bigr|_{Q} := \lim _{t\rightarrow \infty
}\Bigl(\sideset{}{'}\sum_{\substack{x\in L\\|Q (x)|<t}} a (x)  \Bigr),
\end{equation}
where $Q$ is any finite product of real-valued linear forms on $\R
^{n}$ that do not vanish on $\Q ^{n}\smallsetminus \{0 \}$.  As we
explain in \S\ref{qlimit}, this limit depends on $Q$ in a rather simple way.
Nevertheless, to keep notation to a minimum, we assume for now that
the series \eqref{Dede.sum} converges absolutely. 

\subsection{}\label{nature}
The arithmetic nature of the values of $S$ is well known.  Up to
a power of $2\pi i$, they are always rational numbers if $v\in \Q
^{n}$.  However, the explicit calculation of these values is not easy,
especially if $|\det \sigma |$ is large.  In this paper, we 
exhibit a polynomial-time algorithm for calculating $S$ efficiently.
Before stating our main result, we briefly review a few special cases.

In the case $n=1$, $\sigma =1$, $L=\Z$, we have
\begin{equation}
S (L,\sigma ,e,v) = -\frac{(2\pi i)^{e}}{e!}\B _{e} (v),
\end{equation}
where $\B _{e} (v)$ is the periodic function (with period lattice $\Z
$) that coincides with the classical Bernoulli polynomial $B_{e} (v)$
on the interval $0<v<1$.  

More generally, if $L = \Z ^n$, then we
show in Proposition \ref{how.to.sum} that 
\begin{equation}\label{class.dede}
S (L,\sigma ,e,v) = \kappa \sum _{r\in L/\sigma L} \B _{e_{1}} (u_{1})\cdots
\B _{e_{n}} (u_{n}), 
\end{equation}
where 
\begin{equation}
u = \sigma ^{-1} (r+v)\quad \hbox{and}\quad \kappa = \Sign (\det \sigma) \prod
_{j=1}^{n} \frac{(-2\pi i)^{e_{j}}}{e_{j}!}.
\end{equation}
The finite sum on the right of \eqref{class.dede} is a classical
Dedekind sum.  Although theoretically satisfying, from a computational
point of view this sum is of interest only if $|\det \sigma |$ (the
number of summands) is relatively small.

\subsection{}
In general, if $\Rank L <n$, no simple analogue of \eqref{class.dede}
seems to exist except for the special case where $\Rank L$ equals the
number of linear forms $\langle x,\sigma _{j}\rangle $ in
\eqref{Dede.sum} that are not proportional to each other when
restricted to the subspace $L\otimes \R $.  Dedekind sums of this
special type are called \emph{diagonal}.  It is clear that the finite
sum formula \eqref{class.dede} remains valid (modulo obvious modifications)
for all diagonal Dedekind sums.  

Let us further say that a diagonal Dedekind sum is \emph{unimodular}
if the corresponding finite sum has exactly one term.  As a measure of
the deviation from unimodularity, we introduce an integer-valued
function $\|\phantom{S}\|$, called the \emph{index}, on the set of all
Dedekind sums (Definition~\ref{index}).  For example, if $L = \Z ^{n}$, then
$\|S\| = |\det \sigma |$.  In particular, $S$ is said to be unimodular
if and only if $\|S\| = 1$.

\subsection{}
We are now ready to describe our main result.

\begin{theorem}\label{mnthm}
Every Dedekind sum $S (L,\sigma ,e,v)$ can be expressed as a finite
rational linear combination of unimodular diagonal sums.  If 
$n$, $\Rank L$, and $e$ are fixed, then this expression can be
computed in time polynomial in $\log\Index{S}$.  Moreover, the number of
terms in this expression is bounded by a polynomial in $\log\Index{S}$.
\end{theorem}

The special case $n=2$ and $e= (1,1)$ has been known for a very long
time.  This is the case of the classical Dedekind-Rademacher sums,
where the simple form of their reciprocity law combined with the
Euclidean algorithm immediately yields a polynomial-time algorithm for
computing their values (cf. Example~\ref{quadratic}).  Our proof of
Theorem~\ref{mnthm} is similar.  First we prove a general reciprocity
law for the sums $S$.  Applying and specializing this law repeatedly
yields a representation of $S$ as a linear combination of diagonal
sums.  Combining it further with the algorithm of
Ash-Rudolph~\cite{ash.rudolph} finally yields a representation by
unimodular diagonal sums.

\subsection{}
Among the many applications, we review in this paper the problem of
calculating the Eisenstein cocycle~\cite{eisenstein} on $GL_{n} (\Q )$
that, when combined with Theorem~\ref{mnthm}, yields a polynomial-time
algorithm for calculating special values of partial zeta functions of
totally real number fields (\S\S\ref{ec}--\ref{vpzf}).  These special
values are of great interest in view of the many conjectures they are
subject to (Leopoldt conjecture, Brumer-Stark conjecture, etc.).  We
also review the connection between Witten's zeta function and
Dedekind sums, and use our techniques to give some explicit formulas
for special values of these functions (\S\ref{witten}).

\subsection{Acknowledgements}
We thank Gautam Chinta for carefully reading preliminary versions of
this article, and for help with the examples in \S\ref{vpzf}.

%
%

\section{The modular symbol algorithm and Dedekind reciprocity}

\subsection{}
Let $\sigma _{1},\dots ,\sigma _{n}\in \Z ^{n}$
be nonzero primitive points, and let $D=|\det (\sigma _{1},\dots
,\sigma _{n})|$.  Let $w\in \Z ^{n}$ be another nonzero primitive
point, and define 
\[
D_{i} (w) := |\det (\sigma _{1},\dots ,\hat \sigma _{i},\dots ,\sigma
_{n}, w)|,\quad \hbox{$i=1,\dots ,n$.}
\]
The following basic result will play a key role:
\begin{proposition}\label{msa}
\cite{ash.rudolph, barvinok} If $D>1$, then there exists $w\in \Z
^{n}\smallsetminus \{0 \}$ such that
\begin{equation}\label{estimate}
0\leq D_{i} (w) < D^{(n-1)/n},\quad \hbox{$i=1,\dots ,n$,}
\end{equation}
and at least one $D_{i} (w)\not =0$.  Moreover, for fixed $n$ the
point $w$ can be constructed in polynomial time.
\end{proposition}

\begin{proof}
(Sketch) Here we only show that $w$ exists.  Let $P$ be the open
parallelotope
\[
P:=\Bigl\{\sum \lambda _{i}\sigma _{i}\Bigm | |\lambda_{i} |< D^{-1/n}
\Bigr\}. 
\]
Then $P$ is an $n$-dimensional centrally symmetric convex body with
volume $2^{n}$.  By Minkowski's theorem (cf.~\cite[IV.2.6]{fr.tay}),
$P\cap \Z ^{n}$ contains a nonzero point.  This is the desired point $w$.

\end{proof}

\begin{remark}
Ash-Rudolph~\cite{ash.rudolph} show that $w$ satisfies $0\leq D_{i}
(w)<D$.  They also show how to construct $w$ using the Euclidean
algorithm.  The stronger estimate \eqref{estimate} and the statement
about polynomial time are due to Barvinok~\cite[Lemma 5.2]{barvinok}.
In practice, to efficiently construct $w$ such that $D_{i} (w)$ is
small, one may use $LLL$-reduction~\cite[\S3]{experimental}.
\end{remark}

\subsection{}
We now state the Dedekind reciprocity law.  For any nonzero point
$v\in \R ^{n}$, let $v^{\perp }$ be the hyperplane $\{x \mid \langle
v,x\rangle =0 \}$.  Let $Q$ be a finite product of real-valued linear
forms on $\R ^{n}$ that do not vanish on $\Q ^{n}\smallsetminus \{0
\}$ and recall that the notation
\[
\sideset{}{'}\sum _{x\in L} a (x)\Bigr|_{Q}
\]
means the sum is to be evaluated using the $Q$-limit (\S\ref{intro},
\S\ref{qlimit}).

\begin{proposition}\label{recip}
Let $\sigma _{0},\dots ,\sigma _{n}\in \Z ^{n}$ be nonzero.  For
$j=0,\dots ,n$, let $\sigma ^{j}$ be the matrix with columns $\sigma
_{0},\dots ,\hat \sigma _{j},\dots ,\sigma _{n}$.  Fix $L\subseteq \Z
^{n}$, and assume $e = \one := (1,\dots ,1)$.  Then for any $v\in \R
^{n}$, we have the following identity among Dedekind sums:
\begin{equation}\label{recip.law}
\sum _{j=0}^{n} (-1)^{j} S (L,\sigma ^{j},\one )\Bigr|_{Q} =
\sum _{j=0}^{n} (-1)^{j} S (L\cap \sigma _{j}^{\perp },\sigma
^{j},\one )\Bigr|_{Q}. 
\end{equation}

\end{proposition}

\begin{proof}
Let $D_{j}=\det \sigma ^{j}$.
We have an identity for rational functions of $x$
\begin{equation}\label{id1}
\sum _{j=0}^{n}\frac{(-1)^{j}D_{j}}{\prod _{k\not =j}\langle x,\sigma
_{k}\rangle } = 0,
\end{equation}
valid for any $x\in \R ^{n}$ satisfying $\langle x,\sigma _{j}\rangle \not =0$ for
$j=0,\dots ,n$.  To see this, consider the $(n+1)\times (n+1)$ matrix
\begin{equation}\label{mat}
\left(\begin{array}{ccc}
\langle x,\sigma _{0}\rangle &\dots &\langle x,\sigma _{n}\rangle \\
\sigma _{0}&\dots &\sigma _{n}
\end{array} \right).
\end{equation}
This matrix is singular, since the first row is a linear combination
of the others.  Expanding by minors along the top row, and dividing by
$\prod \langle x,\sigma _{k}\rangle $ yields \eqref{id1}.

To pass from \eqref{id1} to \eqref{recip.law}, we need to incorporate
the exponential character and sum over $L$ using $Q$.  There is no
obstruction to doing this, although we must omit terms where
\emph{any} linear form vanishes.  We obtain the expression
\begin{equation}\label{eqqq}
\sum _{j=0}^{n} (-1)^{j}\sum \eee (\langle x,v\rangle
)\frac{D_{j}}{\prod _{k\not =j}\langle x,\sigma _{k}\rangle }\biggr|_{Q} = 0,
\end{equation}
where the inner sum is taken over all $x\in L$ with $\langle x,\sigma
_{k}\rangle \not =0$ for $k=0,\dots ,n$.  

In \eqref{eqqq}, the $j$th inner sum corresponds with the Dedekind sum $S
(L,\sigma ^{j},\one )$ except for the terms with $\langle x,\sigma
_{j}\rangle=0$ and $\langle x,\sigma _{k}\rangle \not =0$ for $k\not
=j$.  In other words, to make the $j$th sum into a Dedekind sum, we
must add
\begin{equation}\label{correct}
\sideset{}{'}\sum_{x\in L\cap \sigma _{j}^{\perp }}\eee(\langle
x,v\rangle )\frac{D_{j}}{\prod _{j\not =k}\langle x,\sigma
_{k}\rangle}\biggr|_{Q}.
\end{equation}
Simultaneously adding and subtracting \eqref{correct} to \eqref{id1}
yields \eqref{recip.law}.
\end{proof}

\subsection{}\label{qlimit}
We now recall the $Q$-limit formula from \cite[Theorem
2]{eisenstein}.  Let 
\begin{equation*}
Q (y) = \prod _{i=1}^{m} Q_{i} (y) 
\end{equation*}
be a product of $m\geq 1$ linear forms
\[
Q_{i} (y) = \sum _{j=1}^{n}Q_{ij}y_{j},
\]
with rationally independent real coefficients $Q_{ij}$.  We think of
$Q$ as an $m\times n$ matrix with rows $Q_{i}$.

Given a vector $e= (e_{1},\dots ,e_{n})$ of positive integers and a
vector $v\in \R ^{n}$, let 
\[
J = \{j\mid \hbox{$e_{j}=1$ and $v_{j}\in \Z$} \}.
\]
If $\#J\equiv 0 \mod 2$, define
\begin{equation}
\BB_{e} (v,Q) = \frac{1}{m}\sum _{i=1}^{m}\left ( \prod _{j\in J}\frac{\Sign Q_{ij}}{2}\right)\prod _{j\not \in J}\B _{e_{j}} (v_{j}),
\end{equation}
otherwise let $\BB _{e} (v,Q) = 0$. In particular, if $J=\varnothing $, then 
\[
\BB_{e} (v,Q) = \prod _{j=1}^{n}\B _{e_{j}} (v_{j}).
\]

Now the $Q$-limit formula can be stated as follows.  Let $\Id_{n}$ be
the $n\times n$ identity matrix.  Then 
\begin{equation}\label{qfmla1}
S (\Z ^{n},\Id_{n},e,v)\Bigr|_{Q} = \kappa \,\BB_{e} (v, Q),
\end{equation}
where
\begin{equation}\label{qfmla2}
\kappa = \prod _{j=1}^{n}\frac{(-2\pi i)^{e_{j}}}{e_{j}!}.
\end{equation}

%
%

\section{Diagonality and unimodularity}

\subsection{}
We begin with some simplifying assumptions to ease the exposition.  We
define the \emph{rank} of $S=S (L,\sigma ,e,v)$ to be the rank of the
lattice $L$.

Suppose that the rank of $S$ is $\ell $, and for any $k$ let $Z
^{k}\subseteq \Z ^{n}$ be the sublattice spanned by the first $k$
standard basis vectors.  We claim that $S$ can be
computed using $Z ^{\ell }$ instead of $L$.  Indeed, writing $L =
gZ^{\ell } $ with a matrix $g\in GL_{n} (\Q ) $, and letting $\sigma '
= g^{t}\sigma $, $v'=g^{t}v$, and $Q' = Qg$, we have
\[
S (L,\sigma ,e,v)\Bigr|_{Q} = (\det g)^{-1} S (Z ^{\ell },\sigma'
,e,v')\Bigr|_{Q'}.
\]
The entries of $\sigma '$ need not be integral, but after multiplying by
an appropriate rational factor we can assume this is true.  In fact,
by further multiplication of $\sigma '$ by rational numbers and
permuting columns, we can write
\[
S (L,\sigma ,e,v)\Bigr|_{Q} = qS (Z ^{\ell },\sigma'
,e,v')\Bigr|_{Q'},\quad q\in \Q ^{\times },
\]
where the pair $(Z ^{\ell },\sigma ')$ satisfies the following:
\begin{enumerate}
\item [(i)] For each column $\sigma_{j} '$, the vector $\sigma _{j}'\cap
Z^{\ell }$ is primitive and integral.
\item [(ii)] If two columns of $\sigma '$ induce proportional linear forms on
$Z ^{\ell }$, then these two linear forms coincide on $Z ^{\ell }$,
and are adjacent columns of $\sigma '$.  
\end{enumerate}

\begin{definition}\label{normalized}
We say that a rank $\ell $ Dedekind sum is \emph{normalized} if the
conditions above are met.
\end{definition}

\subsection{}
Let $S (Z ^{\ell},\sigma ,e,v)$ be a normalized Dedekind sum, and let
$N = \sum e_{j}$.  We claim that without loss of generality, we need
only consider sums for which $e=\one $.  Indeed, let $R ^\ell
\rightarrow R ^{n}\rightarrow \R ^{N}$ be the spans of the first $\ell
$ (respectively $n$) basis vectors in $\R^{N}$, and identify $\R^{n}$
with $R^{n}$.  Let $Q'$ be a product of linear forms on $\R ^{N}$ such
that $Q'$ restricted to $R^{n}$ equals $Q$ restricted to $\R ^{n}$.
We claim that we can construct an $N\times N$ matrix $\sigma '$ and
$v'\in \R ^{N}$ such that
\begin{equation}\label{sum2}
S (Z ^{\ell },\sigma ,e,v)\Bigr|_{Q} = S (Z ^{\ell },\sigma ',\one
,v')\Bigr|_{Q'}.
\end{equation}

To see this, let $\pi \colon \R ^{N} \rightarrow \R ^{\ell }$ be the
projection on the first $\ell $ components.  Given $v\in R ^{n}$, call any $\tilde{v}\in \R ^{N}$
such that $\pi (\tilde{v}) = \pi (v)$ a \emph{lift} of $v$.

Now to construct $S (Z ^{\ell },\sigma ',\one ,v')$, let $v'$ be any
lift of $v$.  If the restriction of the linear form $\sigma _{j}$ to
$Z^{\ell }$ appears on the left of \eqref{sum2} with multiplicity
$e_{j}$, set $j$ columns of $\sigma '$ to be $j$ different lifts of
$\sigma _{j}$.  If we choose these lifts so that $\det \sigma ' = \det
\sigma $, then we obtain \eqref{sum2}.

\begin{definition}\label{prop.embed}
We say that a Dedekind sum with $e=\one $ is \emph{properly
embedded}.
\end{definition}

\subsection{}
Let $S (Z ^{\ell },\sigma ,\one ,v)$ be a properly embedded,
normalized Dedekind sum.  Let $\sets{n}$ be the set $\{1,\dots ,n \}$.

\begin{definition}\label{index}
The \emph{index} of $S$, denoted $\Index{S}$, is defined to be
\[
\Max _{I\subset \sets{n}} \,\bigl|\det \bigl(\pi (\sigma
_{i_{1}}), \dots , \pi (\sigma _{i_{\ell }})\bigr)\bigr|,
\]
where the maximum is taken over all subsets $I = \{i_{1},\dots
,i_{\ell }\}$ of cardinality $\ell $.  A Dedekind sum is called
\emph{unimodular} if $\Index{S} = 1$.
 
\end{definition}

\subsection{}\label{part.sec}
Now define a partition 
\begin{equation}\label{part}
\sets{n} = \bigsqcup_{k=1}^{s} I_{k}, \quad
\hbox{$\ell \leq s\leq n$}
\end{equation}
as follows.  Put
\[
i,j\in I_{k} \quad \hbox{if and only if}\quad \pi
(\sigma _{i}) = \pi (\sigma _{j}).
\]
In other words, two elements of $\sets{n}$ are in the same set of the
partition if the corresponding columns of $\sigma $ induce
the same linear form on $Z ^{\ell }$.

Let $p_{k}=\#I_{k}$.  The vector $p(S) = (p_{1},\dots ,p_{s})$ is
called the \emph{type} of $S$.  To emphasize the type, we 
relabel the columns of $\sigma $ as 
\begin{equation}\label{labelling}
(\sigma _{1}^{1},\dots ,\sigma ^{p_{1}}_{1},\sigma _{2}^{1},\dots
,\sigma _{2}^{p_{2}}, \dots , \sigma _{s}^{1},\dots ,\sigma _{s}^{p_{s}}).
\end{equation}
For any $k=1,\dots ,s$ and any $i=1,\dots, p_{k}$, we denote the point
$\pi (\sigma _{k}^{j})\in \Z ^{\ell }$ by $\sigma _{I_{k}}$.

\begin{definition}\label{diag}
A Dedekind sum is called \emph{diagonal} if $p (S)$ 
has length $\ell$.
\end{definition}

We omit the proof of the following simple lemma:

\begin{lemma}\label{diag.example}
Any normalized, properly embedded rank 1 Dedekind sum is both
diagonal and unimodular.
\end{lemma}

\subsection{}
The virtue of diagonality and unimodularity is the following:

\begin{proposition}\label{how.to.sum}
Keep the notation of the preceding section.
Let $S (Z ^{\ell },\sigma ,\one ,v)$ be properly embedded,
normalized, and diagonal.  Let $\rho$ be the $\ell \times \ell $
matrix $(\sigma _{I_{1}},\dots ,\sigma _{I_{\ell }})$, and let
$(p_{1},\dots ,p_{\ell })$ be the type of $S$.  Then $S|_{Q}$ is
well-defined.  Moreover, 
\begin{equation}\label{diagsum}
S (Z^{\ell },\sigma ,\one ,v)\Bigr|_{Q} = \frac{\kappa \det{\sigma }}{|\det{\rho }|}\sum _{r\in \Z
^{\ell}/\rho \Z^{\ell }} \BB_{p} (u, Q'),
\end{equation}
where
\begin{align}\label{diagsum2}
u &= \rho^{-1}(r+\pi (v)),\quad Q'= Q''\circ \rho^{-t}\circ  \pi, \quad \hbox{and}\\
\kappa &= \prod_{j=1}^{\ell } \frac{(-2\pi i)^{p_{j}}}{p_{j}!}.
\end{align}
Here $Q'' $ is the restriction of $Q$ to $\R^{\ell }$ into $R^{\ell }$. 
If $S$ is unimodular, the sum \eqref{diagsum} has only
one term.
\end{proposition}

\begin{proof}
The second statement follows easily from the first, so we focus on the
first.  By definition, we have
\begin{equation}\label{eqa}
S (Z^{\ell },\sigma ,\one ,v)\Bigr|_{Q} = \sideset{}{'}\sum _{x\in
Z^{\ell }}\eee (\langle x, v\rangle )\frac{\det \sigma }{\prod _{1\leq
k\leq \ell }\langle x, \sigma^{1} _{k}\rangle ^{p_{k}}}\biggr|_{Q}.
\end{equation}
Letting $Q' = Q''\circ \rho^{-t}\circ  \pi$ and $v' = \rho^{-1}\pi (v)$, the right of
\eqref{eqa} becomes 
\begin{equation}\label{eqb}
\det (\sigma )\sideset{}{'}\sum_{y\in \rho^{t}\Z ^{\ell }}\eee (\langle
y,v'\rangle )\frac{\det \rho }{\prod _{1\leq k\leq \ell }
y_{k}^{p_{k}}}\biggr|_{Q'}.
\end{equation}
Inserting the character relations
\begin{equation}\label{eqc}
\sum _{r\in \Z ^{\ell }/\rho\Z ^{\ell }} \eee (\langle y,\rho
^{-1}r\rangle ) = 
\begin{cases}
0,&y\in \Z ^{\ell} \smallsetminus \rho^{t}\Z ^{\ell },\\
\#(\Z ^{\ell }/\rho\Z ^{\ell }),&y\in \rho^{t}\Z ^{\ell },
\end{cases}
\end{equation}
we obtain
\begin{align}\label{eqd}
S (Z^{\ell },\sigma ,\one , v) &=\frac{\det\sigma }{|\det \rho|}\sum _{r\in
\Z ^{\ell }/\rho\Z ^{\ell }}\sideset{}{'}\sum _{y\in \Z ^{\ell }} \frac{\eee (\langle y,\rho^{-1} (r+v)\rangle )}{\prod y_{k}^{p_{k}}}\biggr|_{Q'}\\
&=\frac{\det\sigma }{|\det \rho|}\sum _{r\in\Z ^{\ell }/\rho\Z ^{\ell
}} S (\Z ^{\ell },\Id_{n},p,u)\biggr|_{Q'}\label{eqd:2}
\end{align}
where $u=\rho ^{-1} (r+\pi (v))$.

The proposition follows now from the $Q$-limit formula \eqref{qfmla1}.
\end{proof}

%
%

\section{Algorithms}\label{algorithms}

\subsection{}
In this section we prove that any Dedekind sum is a $\Q $-linear
combination of diagonal, unimodular sums.  We begin with a lemma.  For
simplicity, we abbreviate the Dedekind sum to $S (Z ^{\ell },\sigma
)$.

\begin{lemma}\label{index.no.change}
Let $S (Z ^{\ell },\sigma )$ be a normalized, properly embedded Dedekind sum,
and let $\sigma _{i}$ be a column of $\sigma $.  Then 
\[
\bigl\|S (Z ^\ell\cap \sigma_{i}^{\perp } ,\sigma )\bigr\| \leq
\bigl\|S (Z ^{\ell },\sigma )\bigr\|.
\]
\end{lemma}

\begin{proof}
Without loss of generality, assume that $\sigma _{i} = \sigma _{1}$.
We can represent $\sigma $ as
\begin{equation}\label{mat1}
\sigma =\left( \begin{array}{ccc|ccc|c|cc}
\sigma _{I_{1}}&\cdots&\sigma _{I_{1}}&\sigma _{I_{2}}&\cdots&
\sigma
_{I_{2}}&\cdots &\cdots &\sigma _{I_{s}}\\
*&\cdots&*&*&\cdots&*&\cdots&\cdots&*
\end{array}\right),
\end{equation}
where there are $p_{k}$ columns of the form $(\sigma _{I_{k}},*)^{t}$.
The stars in the last $n-\ell $ rows represent numbers that are
irrelevant, since they don't affect the value of the sum.

Let $\gamma \in GL_{n} (\Q )$ be a matrix that carries $Z ^{\ell }\cap
\sigma _{1}^{\perp }$ onto $Z ^{\ell -1}$.  For $k=1,\dots ,\ell $ let
$\sigma' _{I_{k}}=\gamma \sigma _k^{j}$, where $\sigma _{k}^{j}$ is
any lift of $\sigma _{I_{k}}$.  Then $\bar\gamma \sigma $ has the form
\begin{equation}\label{mat2}
\bar\gamma\sigma = \left( \begin{array}{ccc|ccc|c|cc}
0&\cdots&0&\sigma' _{I_{2}}&\cdots&
\sigma'_{I_{2}}&\cdots &\cdots &\sigma' _{I_{s}}\\
\varepsilon &\cdots&\varepsilon &0&\cdots&0&\cdots&\cdots&0\\
*&\cdots&*&*&\cdots&*&\cdots&\cdots&*
\end{array}\right).
\end{equation}
Here the row blocks have sizes $\ell -1$, $1$, and $n-\ell $, and
$\varepsilon = \pm 1$.

Now let $d$ be the determinant of any $(\ell -1)\times (\ell -1)$ minor
from the top $\ell -1$ rows of \eqref{mat2}.  This
determinant will be the same up to sign as the determinant of an $\ell
\times \ell $ minor containing $\sigma _{1}$ from the top $\ell $ rows
of \eqref{mat1}.  Hence $|d|\leq \Index{S}$,
and the proof is complete.

\end{proof}

Now we come to our first main result.

\begin{theorem}\label{th1}
Let $S = S (Z ^{\ell },\sigma )$ be a normalized, properly embedded
Dedekind sum.  Then we may write
\begin{equation}\label{theq}
S = \sum_{\varrho \in R} q_{\varrho }S (Z^{\ell},\varrho )+\sum_{\tau
\in T} q_{\tau }S (Z^{\ell-1 },\tau),
\end{equation}
where
\begin{enumerate}
\item the sets $R$ and $T$ are finite,
\item $q_{\varrho }, q_{\tau }\in \Q $,
\item each $S (Z^{\ell },\varrho)$ is diagonal, and 
\item each of the Dedekind sums on the right of \eqref{theq} has index
$\leq \Index{S}$.
\end{enumerate}
\end{theorem}

\begin{proof}
Write $\sigma = (\sigma _{1}^{1},\dots ,\sigma ^{p_{1}}_{1},\sigma
_{2}^{1},\dots ,\sigma _{2}^{p_{2}}, \dots , \sigma _{s}^{1},\dots
,\sigma _{s}^{p_{s}})$ as in \eqref{labelling}.  Permuting columns if
necessary, we may assume that $\sigma _{I_{1}},\dots ,\sigma _{I_{\ell
}}$ are linearly independent in $Z ^\ell $.  We will show that we can
write an expression like \eqref{theq} so that the rank $\ell $ sums on the
right have type
\[
(p_{1},\dots ,p_{i}-1,\dots ,p_{\ell },\dots ,p_{s}+1),
\]
for some $1\leq i\leq \ell $.  Iterating this construction proves that
we can have the rank $\ell $ sums on the right of \eqref{theq}
diagonal.

We proceed as follows.  Since $\sigma _{I_{1}},\dots ,\sigma _{I_{\ell
}}$ are linearly independent, we can find unique rational numbers
$\alpha _{j}$ such that
\begin{equation}\label{alphas}
\sigma _{I_{s}} = \sum _{j=1}^{\ell }\alpha _{j}\sigma _{I_{j}}.
\end{equation}

Now we want to apply the relation in Proposition~\ref{recip}.  For
each $\sigma _{I_{j}}$, we choose $p_{j}$ rational lifts
$\tilde{\sigma} ^{1}_{j}$,\dots ,$\tilde{\sigma} _{j}^{p_{j}}$, not
necessarily equal to the columns of $\sigma $, and use them to form a
matrix $\tilde{\sigma }$.  Clearly $S (Z ^{\ell },\tilde{\sigma}) = S
(Z ^{\ell },\sigma )$.  

Now define
\[
\tilde{w} = \sum _{j=1}^{\ell }\alpha _{j}\tilde{\sigma} ^{1}_{j}.
\]
Clearly $\pi (\tilde{w}) = \sigma _{I_{s}}$.  Write
$\tilde{\sigma} _{i}^{j} (\tilde{w})$ for the matrix made from $\sigma $ by
replacing the column $\sigma _{i}^{j}$ with $\tilde{w}$.  
Using the columns of $\tilde{\sigma} $ and $\tilde{w}  $ in
\eqref{recip.law}, we find
\begin{multline}\label{star}
S (Z^{\ell },\tilde{\sigma} ) = S (Z ^{\ell }\cap \tilde{w}^{\perp },\tilde{\sigma} ) + \sum (\varepsilon _{i}^{j}) S
\bigl (Z ^{\ell }, \tilde{\sigma} _{i}^{j} (\tilde{w})\bigr)\\
-\sum (\varepsilon _{i}^{j}) S
\bigl(Z ^{\ell }\cap (\tilde{\sigma} _{i}^{j})^{\perp }, \tilde{\sigma} _{i}^{j} (\tilde{w})\bigr),
\end{multline}
where $\varepsilon _{i}^{j} \in \{\pm 1 \}$ is determined by the
reciprocity law, and the sums are over pairs $(i,j)$ satisfying $1\leq
i\leq s$ and $1\leq j\leq p_{i}$.

We claim that the sums in \eqref{star} actually have only $\ell $
terms.  This follows since the points 
\[
\tilde{\sigma} _{1}^{1},\dots ,\tilde{\sigma} _{\ell }^{1}, \tilde{w}
\]
are dependent.  Hence any sum such that $\tilde{\sigma} _{i}^{j }
(\tilde{w})$ contains these points vanishes.  Moreover, the sum $S (Z
^{\ell }\cap \tilde{w}^{\perp },\tilde{\sigma} )$ is zero, since by
construction a column of $\tilde{\sigma }$ induces a linear form
vanishing on $Z ^{\ell }\cap \tilde{w}^{\perp }$.  Hence \eqref{star}
becomes
\begin{equation}\label{star2}
S (Z^{\ell },\tilde{\sigma} ) = \sum_{i=1}^{\ell } S
\bigl (Z ^{\ell }, \tilde{\sigma} _{i}^{1} (\tilde{w})\bigr)
-\sum_{i=1}^{\ell } S \bigl(Z ^{\ell }\cap (\tilde{\sigma}
_{i}^{1})^{\perp }, \tilde{\sigma} _{i}^{1} (\tilde{w})\bigr). 
\end{equation}

Now consider the types of the rank $\ell $ Dedekind
sums on the right of \eqref{star2}.  If the type of $S$ was
\[
(p_{1},\dots ,p_{i},\dots ,p_{\ell },\dots ,p_{s}),
\]
then the type of $S (\Z ^{\ell },\tilde{\sigma} _{i}^{1} (\tilde{w}))$ is 
\[
(p_{1},\dots ,p_{i}-1,\dots ,p_{\ell },\dots ,p_{s}+1).
\]
Hence by induction we can write $S$ as a finite $\Q $-linear
combination of diagonal rank $\ell$ Dedekind sums plus sums of lower
rank, which proves 1--3 of the statement.

To complete the proof, we must show that the indices of the Dedekind
sums on the right of \eqref{theq} are no larger than $\Index{S}$.
Indeed, Lemma~\ref{index.no.change} implies that the indices on the
right of \eqref{star} are no larger than $\Index{S}$, and so the claim
follows.
\end{proof}

\begin{remark}
There are some simplifications in \eqref{star2} that are worth
mentioning if one wishes to implement the diagonalization algorithm.
First, if we construct $\tilde{\sigma }$ so that $\det \tilde{\sigma
}=1$, then $\det\tilde{\sigma }_{i}^{1} = \alpha _{j}$.  This means
that these determinants can be computed when one computes $\sigma
_{I_{s}}$.  Second, not all the terms on the right of \eqref{star2}
necessarily appear.  In particular, the rank $\ell -1$ sum $S \bigl(Z
^{\ell }\cap (\tilde{\sigma} _{i}^{1})^{\perp }, \tilde{\sigma}
_{i}^{1} (\tilde{w})\bigr)$ appears on the right of \eqref{star2} only
if $p_{i} = 1$.
\end{remark}

\begin{theorem}\label{msa.dede}
With the notation as in Theorem~\ref{th1}, we can write
\begin{equation}\label{lts2}
S = \sum_{\varrho \in R} q_{\varrho }S (Z^{\ell},\varrho )+
\sum_{\tau \in T} q_{\tau }S (Z^{\ell-1 },\tau),
\end{equation}
where the sums of rank $\ell $ on the right are diagonal and
unimodular, and the rank $\ell -1$ sums have index $\leq \Index{S}$.
\end{theorem}

\begin{proof}
By Theorem~\ref{th1}, we may take $S$ to be diagonal.  Suppose that
\[
D = |\det (\sigma _{I_{1}},\dots ,\sigma _{I_{\ell }})| > 1.
\]
Then by
Proposition~\ref{msa}, there exists $w\in Z ^{\ell }$ such that 
\[
D_{j} =|\det (\sigma _{I_{1}},\dots ,\hat \sigma _{I_{j}},\dots
,\sigma _{I_{\ell }}, w)|
\]
satisfies $0\leq D_{j}< D^{(\ell -1)/\ell }$. As in the proof of
Theorem~\ref{th1}, write $w = \sum_{j=1}^{\ell } \alpha _{j}\sigma
_{I_{j}}$, and set $\tilde {w} = \sum_{j=1}^{\ell } \alpha _{j} \sigma
^{1}_{j}$.

Now apply Proposition \ref{recip}, using $\tilde{w}$ and the columns
of $\sigma $, to write $S$ as a finite $\Q $-linear combination of new
Dedekind sums.  These sums won't be diagonal, but we can apply the
proof of Theorem \ref{th1} with $w$ playing the role of $\sigma
_{I_{s}}$.  The resulting sums will include lifts of $w$ in their
columns and will be diagonal.  Thus the resulting rank $\ell $ sums
will have index $< \Index{S}$, and the sums of lower rank will satisfy
the conditions in the statement of Theorem~\ref{th1}.  By induction on
the index, this completes the proof.
\end{proof}

\begin{corollary}\label{cor}
Any Dedekind sum can be written as a finite $\Q $-linear combination
of diagonal, unimodular sums.
\end{corollary}

\begin{proof}
First normalize and embed properly.  By Lemma \ref{diag.example},
any rank 1 Dedekind sum is automatically unimodular and diagonal.  The
result follows by applying Theorems \ref{th1} and~\ref{msa.dede} and
descending induction on the rank.
\end{proof}

%
%

\section{Complexity}\label{complexity}

\subsection{}
In this section we discuss the computational complexity of Corollary
\ref{cor}.  In particular, we show that if $n$ and $\ell $ are fixed,
and $S = S (Z ^{\ell },\sigma )$ is normalized and properly embedded,
then we can form a finite $\Q $-linear combination of diagonal,
unimodular Dedekind sums 
\[
S = \sum _{\substack{\varrho \in R\\
k\leq \ell}} q_{\varrho } S (Z ^{k},\varrho ),
\]
where $\#R$ is bounded by a polynomial in
$\log\Index{S}$.  As a corollary we obtain that this expression can be
computed in polynomial time.  

To do this, we must make a more detailed
analysis of proofs in \S\ref{algorithms}.  We begin by analyzing diagonality.

\begin{lemma}\label{diag.const}
Let $S=S (Z ^{\ell },\sigma )$ be a normalized, properly embedded
Dedekind sum, and write 
\begin{equation}\label{st3}
S = \sum_{\varrho \in R} q_{\varrho }S (Z^{\ell},\varrho )+\sum_{t\in T} q_{\tau }S (Z^{\ell-1 },\tau)
\end{equation}
as in Theorem~\ref{th1}, so that the rank $\ell $ sums in \eqref{st3}
are diagonal.  If $\ell >1$, there exist constants $M_{n,\ell }$ and
$N_{n,\ell }$ such that 
\[
\#R \leq M_{n,\ell } \quad \hbox{and}\quad \#T \leq N_{n,\ell }.
\]
\end{lemma} 

\begin{proof}
Write $p (S) = (p_{1},\dots ,p_{s})$, where $\ell \leq s\leq n$.  By
the proof of Theorem \ref{th1}, we know how to pass from a sum of type 
\begin{equation}\label{typea}
(p_{1},\dots ,p_{i},\dots ,p_{\ell },\dots ,p_{s})
\end{equation}
to a linear combination of sums of types 
\begin{equation}\label{typeb}
(p_{1},\dots ,p_{i}-1,\dots ,p_{\ell },\dots ,p_{s}+1),\quad
\hbox{$i=1,\dots ,\ell$.}
\end{equation}
By iterating this, we pass from $S$ to a linear combination of sums
with types 
\begin{equation}\label{typec}
(p_{1}',\dots ,p'_{i-1},0,p'_{i+1},\dots ,p'_{\ell },\dots ,p'_{s}),
\quad \hbox{$i=1,\dots ,\ell$.}
\end{equation}
We will bound the number of rank $\ell $ (respectively rank $\ell -1$)
sums produced in passing from \eqref{typea} to \eqref{typec} by a
constant $M_{n,\ell }^{(s)}$ (resp.  $N_{n,\ell }^{(s)}$).  We can
then take
\[
M_{n,\ell } = \prod _{s=\ell +1}^{n}M_{n,\ell }^{(s)},
\]
and similarly for $N_{n,\ell }$.

To describe what happens in going from \eqref{typea} to
\eqref{typec}, we use a geometric construction.  Let $B=B (p_{1},\dots
,p_{\ell })$ be the set
\[
B = \bigl\{ (x_{1},\dots ,x_{\ell })\in \Z ^{\ell }\bigm| 0\leq
x_{i}\leq p_{i},  i=1,\dots ,\ell\bigr\}.
\]
The points in $B$ correspond to types of intermediate sums in the
passage from \eqref{typea} to \eqref{typec}.  In particular, passing
from \eqref{typea} to \eqref{typeb} can be encoded by moving from
$(x_{1},\dots ,x_{\ell})$ to $(x_{1},\dots ,x_{i}-1,\dots ,x_{\ell })$
in $B$.  Moreover, the sums of the form \eqref{typec} correspond to
the subset of points $B_{0}\subset B$ with exactly one coordinate $0$.
(See Figure~\ref{b23.fig}.)

\begin{figure}[tbh]
\begin{center}
\includegraphics{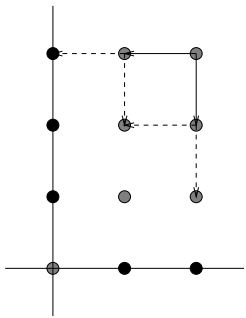}
\end{center}
\caption{\label{b23.fig}The set $B (2,3)$.  The arrows represent
applications of Proposition \ref{recip.law}.  The black dots represent
points in $B_{0}$.}
\end{figure}

Now the constant $M_{n,\ell }^{(s)}$ will be given by $\Max \#B_{0}$,
as $B (p_{1},\dots ,p_{\ell })$ ranges over all possibilities for
fixed $n$, $\ell $, and $s$.  If we allow the $p_{i}$ to become
continuous parameters, then a simple computation shows that the
maximum occurs when $p_{1}=\cdots=p_{\ell } = (n-s+\ell )/\ell $.
With these conditions we have 
\begin{equation}\label{mnls}
M_{n,\ell }^{(s)} \leq \ell \left(\frac{n-s+\ell }{\ell }\right)^{\ell -1}.
\end{equation}

The constant $N_{n,\ell }^{(s)}$ can be computed similarly.  There are
$(\ell +1)$ sums of rank $\ell -1$ produced for each point in 
\[
B_{+} := \bigl\{
(x_{1},\dots ,x_{\ell })\in B \bigm| x_{i} \not = 0\bigr\}.
\]
One finds
again that the maximum occurs when $p_{1}=\cdots=p_{\ell } =
(n-s+\ell )/\ell $, and is
\[
N_{n,\ell }^{(s)} \leq (\ell+1) \left(\frac{n-s+\ell }{\ell }\right)^{\ell}.
\]
\end{proof}

\begin{proposition}\label{polytime1}
If $n$ and $\ell $ are fixed, then \eqref{st3} can be constructed in
constant time, independent of $\Index{S}$.
\end{proposition}

\begin{proof}
This follows easily from the proof of Lemma \ref{diag.const}.  Forming
the expression \eqref{st3} is purely combinatorial, and makes no
reference to $\Index{S}$.  In particular, the number of steps needed
can be bounded for fixed $n$ and $\ell $.
\end{proof}

\subsection{}
Now we investigate the size of the output in Theorem \ref{msa.dede}.
The main step of the proof of Theorem \ref{msa.dede} shows how given
$S$, one may write 
\begin{equation}\label{lts}
S = \sum_{\varrho \in R} q_{\varrho }S (Z^{\ell},\varrho )+ \sum_{\tau
\in T} q_{\tau }S (Z^{\ell-1 },\tau),
\end{equation}
where the sums of rank $\ell $ on the right satisfy $\Index{S
(Z^{\ell},\varrho)} < \Index{S}^{(n-1)/n}$.  Denote by $C_{\ell }$
(respectively $C_{\ell -1}$) the number of rank $\ell $ (resp. rank
$\ell -1$) sums on the right of \eqref{lts}.

\begin{lemma}\label{computingC}
We have $C_{\ell } \leq M^{(\ell +1)}_{n+1,\ell }$ and $C_{\ell -1}
\leq N_{n+1,\ell }^{(\ell +1)}. $
\end{lemma}

\begin{proof}
The proof is very similar to the that of Lemma \ref{diag.const}, with
the following twist: we begin with a diagonal sum, increase the number
of distinct linear forms by one, and then make the sums diagonal
again.  We can keep track of the number of sums produced using the set
$B$ as in the preceding proof, although we must replace $n$ with $n+1$
to accommodate the extra initial step.  We leave the details to the
reader.
\end{proof}

\subsection{}
Now consider the expression \eqref{lts}.  To complete the proof of
Theorem \ref{msa.dede}, we repeat the process that produced
\eqref{lts} until all rank $\ell $ sums are diagonal and unimodular,
and we obtain the expression \eqref{lts2}.   Using Lemma
\ref{computingC} and the estimate \eqref{estimate}, we can bound the
number of rank $\ell $ sums produced.

\begin{proposition}\label{poly1}
(cf. \cite[Theorem 5.4]{barvinok}) Write $S$ as a sum of diagonal,
unimodular rank $\ell $ sums and lower rank sums as in \eqref{lts2}.
Then the number of rank $\ell $ sums on the right of \eqref{lts2} is
bounded by
\begin{equation}\label{bound1}
C' (\log\Index{S})^{A\log C_{\ell }},
\end{equation}
where $A = \log (n/ (n-1))^{-1}$, and $C'$ is a constant independent of
$S$.
\end{proposition}

\begin{proof}
In \eqref{lts}, we have 
\[
0\leq \Index{S (Z ^{\ell },\varrho )} < \Index{S}^{(n-1)/n}.
\]
Thus after $t$ iterations we'll have 
\[
0\leq \Index{S (Z ^{\ell },\varrho )} < \Index{S}^{( (n-1)/n)^{t}}.
\]
Since the index of a Dedekind sum is always an integer, the condition
for termination is that for some $\varepsilon >0$, we have
\begin{equation}\label{test}
\Index{S}^{( (n-1)/n)^{t}} \leq 2-\varepsilon, \quad \hbox{or} \quad
t\geq \frac{\log\log \Index{S}-\log\log (2-\varepsilon )}{\log n - \log (n-1)}.
\end{equation}

On the other hand, by Lemma \ref{computingC}, we know that $t$
iterations will produce no more than $C_{\ell }^{t}$ sums of rank
$\ell $.  So fix $\varepsilon >0$, set $A = \log n - \log (n-1)$, and
let 
\[
C'_{\ell } = \exp\left(\!\left (\frac{-\log\log (2-\varepsilon
)}{A}+1\right)\log C_{\ell }\right).
\]
Then 
\[
C_{\ell }^{t}\leq C'_{\ell } \bigl(\log \Index{S}\bigr)^{A\log C_{\ell }}.
\]
Now if we define $C' = \Max _{\ell \leq n} C_{\ell }'$, we obtain
\eqref{bound1}. 
\end{proof}

\subsection{}
We are now ready to discuss the complexity of our
algorithms. 

\begin{theorem}\label{th3}
Let $S=S (Z ^{\ell },\sigma )$ be a normalized, properly embedded
Dedekind sum.  Using Corollary~\ref{cor}, write $S$ as a $\Q $-linear
combination of diagonal unimodular sums.  Then there exists a
polynomial $P_{n,\ell }$ such that the number of terms in the
expression is bounded by $P_{n,\ell } (\log \Index{S})$.  Moreover, we
have
\[
\deg P_{n,\ell } \leq A \log\left (\frac{\ell !\,n^{\ell (\ell
-1)/2}}{2^{1}3^{2}\cdots\ell^{(\ell -1)}}\right),
\]
where $A = \log (n/ (n-1))^{-1}$.
\end{theorem}

\begin{proof}
Fix $n$. We proceed by induction on $\ell $.

First, if $\ell =1$, then by Lemma \ref{diag.example} the sum $S (Z^{1}
,\sigma )$ is already diagonal and unimodular.  Hence we may take
$P_{n,1} \equiv 1$.

Next, assume that the statement is true for sums of rank $\ell -1$, and
let $P_{n,\ell -1}$ be the corresponding polynomial.  First we claim
that without loss of generality, we need only consider the case that
$S$ is diagonal.  Indeed, apply 
Theorem~\ref{th1} and write
\begin{equation}\label{st4}
S = \sum_{\varrho \in R} q_{\varrho }S (Z^{\ell},\varrho )+\sum_{\tau
\in T} q_{\tau }S (Z^{\ell-1 },\tau),
\end{equation}
where the rank $\ell $ Dedekind sums are diagonal, and all the
Dedekind sums have index $\leq \Index{S}$.  
By Lemma \ref{diag.const}, the sets $R$ and $T$ have a bounded number
of elements independent of $S$.  Hence we
may bound the output for diagonal $S$, and then multiply this by a
constant to obtain our final answer.

Now we apply Theorem~\ref{msa.dede}, and we must count the number of
rank Dedekind sums produced.  By Proposition \ref{poly1}, we know
that the total number of rank $\ell $ sums will be bounded by 
\begin{equation}\label{est1}
C' (\log\Index{S})^{A\log C_{\ell }}.
\end{equation}
Furthermore, each sum of lower rank produced in the proof of Theorem
\ref{msa.dede} can be written as a sum of $\leq P_{n,\ell -1}
(\Index{S})$ Dedekind sums by induction.  To find the total output, we
must count these lower rank sums.

We can do this as follows.  Let $Q = Q_{n,\ell }:= C_{\ell
-1}P_{n,\ell -1}$, where $C_{\ell -1}$ is the constant in Lemma
\ref{computingC}.  Represent the process of reducing the diagonal sum
$S$ to unimodularity by the following diagram:
\begin{equation}\label{diagram2}
\xymatrix{{\hbox{rank $\ell $:}}&1\ar[dr]\ar[r]&C_{\ell }\ar[dr]\ar[r]&C_{\ell }^{2}\ar[dr]\ar[r]&\cdots\ar[dr]\ar[r]&C_{\ell }^{t}\\
{\hbox{rank $\ell-1 $:}}&&Q&QC_{\ell }&\cdots&QC_{\ell }^{t-1}}
\end{equation}
The top row represents the bound on the number of rank $\ell $ sums at
each step of the algorithm, and the bottom
row is the number of rank $\ell -1$ sums.  

According to \eqref{diagram2}, we have produced 
\[
Q\sum _{i=0}^{t-1}C_{\ell }^{i} = Q\,\frac{C_{\ell }^{t}-1}{C_{\ell
}-1} \leq Q\,\frac{C'(\log\Index{S})^{A\log C_{\ell }}-1}{C_{\ell
}-1}
\]
sums of ranks $<\ell $.  Adding this estimate to \eqref{est1}, we
find that the total number of Dedekind
sums produced will be a polynomial of
degree 
\[
\deg Q + A\log C_{\ell } = \deg P_{n,\ell-1 } + A\log C_{\ell }.
\]

To complete the proof, we compute the degree of $P_{n,\ell }$ by
induction.  Indeed, using the estimate
\[
C_{\ell } \leq M_{n+1,\ell }^{(\ell +1)} = \ell\left(\frac{n}{\ell }\right)^{\ell -1},
\]
and that $\deg P_{n,1} = 0$,
an easy computation shows that 
\[
\deg P_{n,\ell } \leq A \log\left (\frac{\ell !\,n^{\ell (\ell
-1)/2}}{2^{1}3^{2}\cdots\ell^{(\ell -1)}}\right),
\]
as required.
\end{proof}

\begin{example}\label{table}
Here is a table of the bound of $\deg P_{n,\ell }$ for small values of $n$
and $\ell $. 
\vskip2em
\begin{center}
\begin{tabular}{l||ccccccccc}
$n,\ell $&2&3&4&5&6&7&8\\
\hline\hline
2 & 1 &&&&&&\\
3 & 2 & 5 &&&&&&\\
4 & 4 & 10 & 15 &&&&\\
5 & 7 & 16 & 25 & 33 &&&\\
6 & 9 & 23 & 37  & 50  & 60 &&\\
7 & 12 & 30 & 50 & 69 & 86 & 99&\\
8 & 15 & 38 & 64 & 90 & 114 & 135 & 150
\end{tabular}
\end{center}
\vskip2em
\end{example}

\begin{corollary}\label{polytime2}
Keeping the same notation as in Theorem \ref{th3}, for fixed $n$ and
$\ell $ we may express $S$ as a finite $\Q $-linear combination of
diagonal unimodular sums in time polynomial in $\log\Index{S}$.
\end{corollary}

\begin{proof}
First, the vector $w$ constructed in Proposition \ref{msa} can be found in
polynomial time in the size of the coefficients of $S$.  In fact,
investigation of \cite[Lemma 5.2]{barvinok} shows that the rational
numbers $\alpha _{j}$ in \eqref{alphas} in the proof of Theorem
\ref{th1} can also be constructed in polynomial time.  This implies that
$w$ and the $\alpha _{j}$ can be found in time polynomial in
$\log\Index{S}$.  The proof of Theorem \ref{th3} then shows that the
final expression can be computed in polynomial time.
\end{proof}

%
%

\section{The Eisenstein cocycle}\label{ec}

\subsection{}
In this section, we briefly review the construction of the Eisenstein
cocycle introduced in \cite{eisenstein}.  In particular, we show that
this is a finite object that can be calculated effectively using
Corollary \ref{cor}.  Roughly speaking, the Eisenstein cocycle
represents a generalization of the classical Bernoulli polynomial
within the arithmetic of the unimodular group $\Gamma =GL_{n}(\Z) $.

\subsection{}
Let $\A = (A_{1},\dots ,A_{n})$ be an $n$-tuple of matrices $A_{i}\in
GL_{n} (\R)$.  For an $n$-tuple $d = (d_{1},\dots ,d_{n})$ of integers
$1\leq d_{i}\leq n$, let $\A (d)\subseteq \R ^{n}$ be the subspace
generated by all columns $A_{ij}$ such that $j<d_{i}$.  (Here $A_{ij}$
denotes the $j$th column of $A_{i}$.)  Writing $\A (d)^{\perp}$ for
the orthogonal complement of $\A (d)$ in $\R^{n}$, we let 
\begin{equation}\label{one}
X (d) = \A (d)^{\perp } \smallsetminus \bigcup_{i=1}^{n} \sigma
_{i}^{\perp }, \quad \hbox{where $\sigma _{i} = A_{id_{i}}$.}
\end{equation}
The $n$-tuple $\A $ determines then the stratification 
\begin{equation}
\R ^{n} \smallsetminus \{0 \} = \bigsqcup_{d\in D} X (d),
\end{equation}
indexed by the finite set 
\[
D = D (\A ) = \{d \mid X (d) \not = \varnothing  \}.
\]
Associated to this decomposition is the rational function $\psi (\A )$
on $\R ^{n}\smallsetminus \{0 \}$ defined by 
\[
\psi (\A ) (x) = \frac{\det (\sigma _{1},\dots ,\sigma _{n})}{\langle
x,\sigma _{1}\rangle \cdots \langle x,\sigma _{n}\rangle },\quad
\hbox{if $x\in X (d)$.}
\]

\subsection{}
More generally, if $P (X_{1},\dots ,X_{n})$ is any homogeneous
polynomial, we form the differential operator $P (-\partial
_{x_{1}},\dots ,-\partial _{x_{n}})$ in the partial derivatives
$\partial _{x_{i}} := \partial /\partial x_{i}$, and define 
\[
\psi (\A ) (P,x) = P (-\partial
_{x_{1}},\dots ,-\partial _{x_{n}}) \psi (\A ) (x).
\]

The last expression can be written more explicitly as 
\begin{equation}
\psi (\A ) (P,x) = \det (\sigma )\sum _{r}P_{r} (\sigma )\prod
_{j=1}^{n} \frac{1}{\langle x,\sigma _{j}\rangle ^{1+r_{j}}},
\end{equation}
where $r$ runs over all decompositions of $\deg (P) =
r_{1}+\cdots+r_{n}$ into nonnegative parts $r_{j}\geq 0$, and $P_{r}
(\sigma )$ is the homogeneous polynomial in the $\sigma _{ij}$ defined
by the expansion 
\[
P (X\sigma ^{t}) = \sum _{r}P_{r} (\sigma ) \prod
_{j=1}^{n}\frac{X_{j}^{r_{j}}}{r_{j}!}.
\]
In the excluded case $x = 0$, it is convenient to set $\psi (\A )
(P,0) = 0$.

The definition of the Eisenstein cocycle $\Psi $ is now easy to state:
\begin{equation}\label{four}
\Psi (\A ) (P,Q,v) := (2\pi i)^{-n-\deg P}\sum _{x\in \Z ^{n}} \eee (\langle x,v\rangle )
\psi (\A ) (P,x)\Bigr|_{Q}.
\end{equation}

The series on right converges provided all components $A_{i}$ of $\A $
are in $GL_{n} (\Q )$.  However, since the convergence is only
conditional, we are forced to introduce the additional parameter $Q$
specifying the limiting process.

\subsection{}
Let $M$ be the set of all complex-valued functions $f (P,Q,v)$ with
$P$, $Q$, $v$ as above ($v\in \R ^{n}$).  Then $M$ is a left $\Gamma
$-module under the action
\[
Af (P,Q,v) = \det ( A) f (A^{t}P, A^{-1}Q, A^{-1}v), \quad \hbox{$A\in
\Gamma$,} 
\]
where the implied $\Gamma $-action on homogeneous polynomials is given
by $(AP) (X) = P (XA)$.  With respect to this action, the map 
$\Psi \colon \Gamma ^{n}\rightarrow M$ has the property
\begin{align}
\Psi &(A\A ) = A\Psi (\A ), \quad\hbox{$A\in \Gamma $, $\A \in \Gamma ^{n}$,}\\
&\sum _{i=0}^{n} \Psi (A_{0},\dots ,\hat A_{i},\dots ,A_{n})=0, \quad
\hbox{$A_{i} \in \Gamma $.}
\end{align}
In other words, $\Psi $ is a homogeneous cocycle on $\Gamma $.  It is
known that $\Psi $ represents a nontrivial cohomology class in
$H^{n-1} (\Gamma ; M)$ \cite[Theorem 4]{eisenstein}.

Combining \eqref{one}--\eqref{four}, we see that $\Psi $ is a finite
linear combination of Dedekind sums, 
\[
\Psi (\A ) (P,Q,v) = (2\pi i)^{-n-\deg P}\sum _{d\in D}\sum _{r} P_{r}
(\sigma )S (L,\sigma ,e,v)\Bigr|_{Q}.
\]
Here $\sigma $ is the matrix with columns $A_{id_{i}}$ for $i=1,\dots ,n$,
$L$ is the lattice $\A (d)^\perp \cap \Z ^{n}$, and $e_{j} = 1+r_{j}$.
The case $P=1$ is of special interest:
\[
\Psi (\A ) (1,Q,v) = (2\pi i)^{-n}\sum _{d\in D} S (L,\sigma ,\one
,v)\Bigr|_{Q}.
\]  
This case yields the classical Dedekind-Rademacher sums if $n=2$, and,
more importantly, it corresponds to special values of partial zeta
functions at $s=0$.

%
%

\section{Values of partial zeta functions}\label{vpzf}

\subsection{}
Let $F$ be a totally real number field of degree $n$ over $\Q $, and
let $\fff$, $\bbb$ be two relatively prime ideals in the ring of
integers $\O _{F}$.  The partial zeta function to the ray class $\bbb
\mod \fff$ is defined by 
\[
\zeta (\bbb, \fff, s) := \sum _{\aaa \equiv \bbb \mod \fff } N (\aaa
)^{-s}, \quad \hbox{$\Re (s)>1$,}
\] 
where $\aaa $ runs over all all integral ideals in $\O _{F}$ such that
the fractional ideal $\aaa \bbb ^{-1}$ is a principal ideal generated
by a totally positive number in the coset $1+\fff \bbb ^{-1}$.
According to Klingen-Siegel, the special values
$\zeta (\bbb ,\fff ,1-s)$, where $s=1,2,3,\dots, $ are well-defined
rational numbers.  In this section, we give a formula for calculating
these numbers in terms of the Eisenstein cocycle $\Psi $.

\subsection{}
The formula depends on the choice of a $\Z $-basis $W$ for the
fractional ideal $\fff \bbb ^{-1} = \sum \Z W_{j}$, together with the
dual basis $W^{*}$ determined by $\Trace (W_{i}^{*}W_{j}) = \delta
_{ij}$.  Here we identify $\alpha \in F$ with the row vector $(\alpha
^{(1)}, \dots , \alpha ^{(n)})\in \R ^{n}$, where the $\alpha ^{(j)}$
are the $n$ different embeddings of $\alpha $ into the field of real
numbers.  Then $W$ can be identified with a matrix in $GL_{n} (\R )$
whose $j$th row is the basis vector $W_{j}$.  Let 
\begin{align*}
P (X)&=N (\bbb )\prod _{i}\sum _{j}X_{j}W_{j}^{(i)},\\
Q (X)&= \prod _{i}\sum _{j}X_{j}(W_{j}^{*})^{(i)},
\end{align*}
and let $v\in \Q ^{n}$ be defined by $v_{j} =
\Trace (W^{*}_{j})$.

The formula also depends on the choice of generators $\varepsilon
_{1},\dots ,\varepsilon _{\nu }$, where $\nu = n-1$, for the group
$U\subset \O _{F}^{\times }$ of totally positive units.  Using the regular
representation $\rho \colon U\rightarrow \Gamma $, defined via $\rho
(\varepsilon ) = W\delta (\varepsilon )W^{-1}$, where $\delta
(\varepsilon )$ is the matrix $\Diag (\varepsilon ^{(1)},\dots
,\varepsilon ^{(n)})$, we identify the units $\varepsilon _{j}$ with
elements $A_{j} = \rho (\varepsilon_{j} )^{t} \in \Gamma $.  (Note
that $\rho $ is the \emph{row} regular representation.)

\subsection{}
Using the bar notation
\[
[A_{1}|\cdots | A_\nu ] := (1,A_{1},A_{1}A_{2},\dots ,A_{1}\cdots A_{\nu })\in \Gamma ^{n},
\]
we have the following proposition expressing the zeta values in terms
of the Eisenstein cocycle:
\begin{proposition}\label{zeta-values}
Let $U_{\fff }$ be the subgroup $U\cap (1+\fff )$, and let $\pi $ run
through all permutations of \hbox{$\{1,\dots ,\nu \}$}.  Then for
$s=1,2,3,\dots$,
\[
\zeta (\bbb ,\fff ,1-s) = \eta \sum _{\varepsilon \in U/U_{\fff }}
\sum _{\pi } \Sign (\pi) \Psi ([A_{\pi (1) }|\cdots | A_{\pi (\nu )}])
(P^{s-1},Q,\rho (\varepsilon )^{t}v).
\]
Here the sign $\eta = \pm
1$ is determined by 
\[
\eta = (-1)^{\nu }\Sign (\det W) \Sign (R), 
\]
where $R = \det (\log \varepsilon _{i}^{(j)})$, $1\leq i,j\leq \nu $.
\end{proposition}

\begin{proof}
This follows from \cite[Corollary, p. 595]{eisenstein} by writing the
fundamental cycle of $U_{\fff }$ in terms of the $A_{j}$.
\end{proof}

\begin{example}\label{quadratic}
We work out the above formula in the case of a real quadratic field
$F$.  Let $\varepsilon >1$ be the fundamental unit of $U$, the group of
totally positive units in $F$, and let $\Z w_{1} + \Z w_{2} = \fff
\bbb ^{-1}$ be a $\Z $-basis of $\fff \bbb ^{-1}$.  Such a basis
determines a matrix $A=\left(\begin{smallmatrix}
a&b\\
c&d
\end{smallmatrix} \right)\in SL_{2} (\Z )$ and a vector $v\in \Q ^{2}$ via 
\[
\left(\begin{array}{c}
\varepsilon w_{1}\\
\varepsilon w_{2}
\end{array} \right) = 
\left(\begin{array}{cc}
a&c\\
b&d
\end{array} \right)
\left(\begin{array}{c}
w_{1}\\
w_{2}
\end{array} \right), \quad v_{1}w_{1} + v_{2}w_{2}= 1.
\]
In addition, we get the normforms
\[
P (X) = N (x_{1}w_{1} + x_{2}w_{2}), \quad Q (X) = N (x_{1}w_{2}-x_{2}w_{1}).
\]
Let $p$ be the smallest positive integer such that $(A^{p}-1)v\in \Z
^{2}$.  Then 
\begin{equation}\label{eqq1}
\zeta_{F} (\bbb ,\fff ,1-s) = \eta \sum _{k\mod p} \Psi (1,A) (P^{s-1},Q,A^{k}v),
\end{equation}
where $\eta = \Sign (w_{2}w_{1}^{(1)}- w_{1}w_{2}^{(1)})$, and, if $s=1$,
\begin{align}\label{}
\Psi \left(\left(\begin{array}{cc}
1&0\\
0&1
\end{array} \right), \left(\begin{array}{cc}
a&b\\
c&d
\end{array} \right) \right) (1,Q,v) &= \frac{a}{2c}\B _{2} (v_{2}) +
\frac{d}{2c}\B _{2} (cv_{1}- av_{2}) \\
&- \sum _{j\mod c} \B _{1}
(\frac{j+v_{2}}{|c|})\B _{1} (a\frac{j+v_{2}}{c}-v_{1})\label{classicalsum},
\end{align}
with an additional correction term $-\Sign (c)/4$ on the right if $v\in \Z ^{2}$.
The finite sum \eqref{classicalsum} is the classical Dedekind-Rademacher sum
$(2\pi i)^{-2}S \left(\Z ^{2}, \left(\begin{smallmatrix}
1&a\\
0&c
\end{smallmatrix} \right), \one ,v\right)\bigr|_{Q}$.  Note that the
number of terms in that sum equals $|c| = | (\varepsilon
-\varepsilon ')/(w-w')|$, where $w = w_{2}/w_{1}$, and the prime is
Galois conjugation.  Depending on
$\varepsilon $, this number can be very large.  To get a more
efficient formula for calculating $\Psi $, we apply the euclidean
algorithm to the first column of $A$ an obtain a product decomposition
\[
A = B_{1}\cdots B_{t}, \quad t\geq 1, \quad B_{j} = \left(\begin{array}{cc}
b_{j}&-1\\
1&0
\end{array} \right).
\]
Then
\begin{equation}\label{eqq2}
\Psi (1,A) = \sum _{\ell =0}^{t-1} (B_{1}\cdots B_{\ell })\Psi
(1,B_{\ell +1}).  
\end{equation}
Here $t$ is roughly $\log|c|$.  Thus the number of terms is
effectively reduced from $|c|$ to $\log|c|$, since
\begin{multline}\label{eqq3}
\Psi \left(\left(\begin{array}{cc}
1&0\\
0&1
\end{array} \right), \left(\begin{array}{cc}
b&-1\\
1&0
\end{array} \right) \right) (P^{s-1}, Q, v) \\
= 
\sum _{r}\left[bP_{r}\left(\begin{array}{cc}
0&b\\
1&1
\end{array} \right)\frac{\B _{2s} (v_{2})}{(2s)!} + P_{r}\left(\begin{array}{cc}
1&b\\
0&1
\end{array} \right)\frac{\B _{1+r_{1}} (v_{1}-bv_{2})}{(1+r_{1})!}\frac{\B _{1+r_{2}} (v_{2})}{(1+r_{2})!}\right],
\end{multline}
where the rational numbers $P_{r}\bigl (\begin{smallmatrix}\alpha &\beta \\
\gamma &\delta \end{smallmatrix}\bigr)$ are the coefficients of the
polynomial
\[
P (\alpha X_{1}+\beta X_{2}, \gamma X_{1}+\delta X_{2})^{s-1} = \sum_{r}P_{r}\left(\begin{array}{cc}
\alpha &\beta \\
\gamma &\delta 
\end{array} \right)\frac{X_{1}^{r_{1}}}{r_{1}!}\frac{X_{2}^{r_{2}}}{r_{2}!}.
\]
In the exceptional case $s=1$, $\fff = (1)$, the correction term
\[
-\frac{1}{8}\{\Sign (w+b) + \Sign (w'+b) \}
\]
must be added to the right side of \eqref{eqq3}.

As a numerical example, we choose $F=\Q (\sqrt{5})$, $\varepsilon =
(3+\sqrt{5})/2$, $\fff = \bbb = (1)$, $w_{1}=-\varepsilon $, $w_{2}=1$.  Then
$\eta = +1$, $A = \left(\begin{smallmatrix}3&-1\\
1&0\end{smallmatrix} \right)$, $P_{11}\left(\begin{smallmatrix}1&3\\
0&1\end{smallmatrix} \right) = 3$, while $P_{20}\left(\begin{smallmatrix}0&3\\
1&1\end{smallmatrix} \right) = 2$, $P_{11}\left(\begin{smallmatrix}0&3\\
1&1\end{smallmatrix} \right) = -7$, $P_{02}\left(\begin{smallmatrix}0&3\\
1&1\end{smallmatrix} \right) = 2$.  Hence, according to \eqref{eqq1}
and \eqref{eqq3}, we get for the value of the Dedekind zeta function
of $F$ at $s=-1$,
\[
\zeta _{F} (-1) = \zeta ((1),(1),-1) = 3 (2-7+2) (-\frac{1}{720})
+ 3 (\frac{1}{12})^{2} = \frac{1}{30}.
\]
\end{example}

\begin{example}\label{cubic}
As a second example, we consider the cubic field $\Q (\theta )$ of
discriminant $148$ given by $\theta ^{3} - \theta ^{2} -3\theta +1=0$.
According to \cite{khan}, the group of totally positive units $U$ is
generated by $\varepsilon _{1} = -3\theta ^{2}+2\theta +10$ and
$\varepsilon _{2} = 5\theta ^{2}+6\theta -2$.

Let $\fff = (2)$ and $\bbb = (1)$.
Then $\fff \bbb ^{-1} = \Z w_{1}+ \Z w_{2} + \Z w_{3}$, where $w_{1} =
2$, $w_{2} = 2 \theta $, and $w_{3} = 2\theta ^{2}$.
With respect to this basis, we find
\[
A_{1} = \left(\begin{array}{ccc}
{10}&{3}&{1}\\
 {2}&{1}&{0}\\
 {-3}&{-1}&{0}
\end{array} \right),
A_{2} = \left(\begin{array}{ccc}
 {-2}&{-5}&{-11}\\
 {6}&{13}&{28}\\
 {5}&{11}&{24}
\end{array} \right),
A_{1}A_{2} = \left(\begin{array}{ccc}
 {3}&{0}&{-2}\\
 {2}&{3}&{6}\\
 {0}&{2}&{5}
\end{array} \right).
\]
Then $v = (1/2,0,0)^{t}$ and $\eta = 1$.  Let $V\subset \Q ^{3}$ be a
complete set of representatives for the orbit of $v+\Z^3$ under the
action of $U$ (via $A_{1}$ and $A_{2}$) on $\Q^3/\Z^3$. Note that $V$
is a finite set.  Thus
\begin{equation}\label{eqq4}
\zeta (\bbb , \fff , 0) = \sum_{v\in V} \bigl (\Psi(1, A_{1},
A_{1}A_{2}) - \Psi (1, A_{2}, A_{1}A_{2})\bigr) (1, Q, v).
\end{equation}

Each term on the right of \eqref{eqq4} breaks up into $10$ Dedekind
sums: one of rank $3$, three of rank $2$, and six of rank $1$.  Note
that the rank $3$ and rank $1$ sums are diagonal, whereas the rank $2$
sums are not.  After making all sums diagonal, we find that each $\Psi
$ in \eqref{eqq4} has $30$ terms.  Applying the summation formula for
diagonal sums (Proposition \ref{how.to.sum}), we see that to evaluate
$\Psi $ on any element of $V$, we must sum $76$ terms.

Since $\varepsilon _{2}^{2} \equiv \varepsilon _{1}\varepsilon
_{2}\equiv 1 \mod \fff $, we can take $V = \{v,  A_{2}v\}$.
This yields $152$ terms altogether, all of which sum to
$0$, and thus $\zeta ((1), (2), 0) = 0$.  This agrees with
\cite{khan}, and also with the observation that since $-1$ preserves
the congruence class of $1+\fff $, and the norm of $-1$ is $-1$, the
special value at $s=0$ must vanish.
 
Now let $\fff = (3)$.  Then we may take $A_{1}, A_{2}$ as above, and
$v = (1/3,0,0)^{t}$ and $\eta =1$.  Since $\varepsilon _{2}^{13}
\equiv \varepsilon _{1}\varepsilon _{2}^{5}\equiv 1 \mod \fff $, we
must sum $13\cdot 76= 988$ terms, and we find $\zeta ((1), (3), 0) =
2/3$, again in agreement with \cite{khan}.

Note that to compute $\zeta ((1), (N), 0)$ for various $N\in \Z $, we
must only compute $A_{1}$, $A_{2}$, and thus $\Psi
(1,A_{1},A_{1}A_{2}) - \Psi (1,A_{2},A_{1}A_{2})$ once.  After this it
is routine to compute special values at $s=0$, and the complexity in
\eqref{eqq4} comes from $\#V$, which can be large, even for small
values of $N$.  A table of $\zeta ((1), (N), 0)$ for several rational
integers $N$ is given below.  The values for $N=2,3,5,7$ are also in
\cite{khan}.

\vskip2em
\begin{center}
\begin{tabular}{|c|r||c|r||c|r||c|r||c|r|}
\hline
$N$&$N\cdot \zeta$&$N$&$N\cdot \zeta$&$N$&$N\cdot \zeta$&$N$&$N\cdot \zeta$&$N$&$N\cdot \zeta$\\
\hline\hline
&&$11$&$-18$&$21$&$-78$&$31$&$74$&$41$&$382$\\
$2$&$0$&$12$&$5$&$22$&$68$&$32$&$15$&$42$&$-228$\\
$3$&$2$&$13$&$-22$&$23$&$12$&$33$&$62$&$43$&$-366$\\
$4$&$1$&$14$&$-20$&$24$&$23$&$34$&$50$&$44$&$1$\\
$5$&$-4$&$15$&$42$&$25$&$106$&$35$&$-54$&$45$&$254$\\
$6$&$4$&$16$&$7$&$26$&$-24$&$36$&$-43$&$46$&$6$\\
$7$&$2$&$17$&$100$&$27$&$-190$&$37$&$20$&$47$&$-570$\\
$8$&$3$&$18$&$-32$&$28$&$25$&$38$&$22$&$48$&$-13$\\
$9$&$-10$&$19$&$82$&$29$&$242$&$39$&$156$&$49$&$-222$\\
$10$&$-2$&$20$&$4$&$30$&$6$&$40$&$2$&$50$&$178$\\
\hline
\end{tabular}
\end{center}
\begin{center}
$\zeta ((1),(N),0)$ for the cubic field of discriminant 148.
\end{center}
\vskip2em
\end{example}

%
%

\section{Witten's zeta function}\label{witten}

\subsection{}
To give another illustration, we recall the definition of Witten's zeta
function \cite{witten}, and show how our algorithms can be used to
compute special values of this function at even integers.  For
unexplained notions from representation theory, the reader may consult
\cite{fulton.harris}.

Let $\ggg$ be a simple complex lie algebra, and let $R$ be the
associated root system.  Let $R^{+}$ (respectively $R^{-}$) be a
subset of positive roots (resp. negative roots), and let $\Delta
\subset R^{+}$ be the set of simple roots.

The roots $R$ generate a lattice $\Lambda _{R}$ in an $\ell
$-dimensional real vector space $E$ endowed with an inner product
$(\phantom{a},\phantom{a})$.  Let $\Lambda _{W}$ be the weight
lattice, which is the dual of $\Lambda _{R}$ with respect to this
inner product.  We denote by $\Omega \subset E$ the set of fundamental
weights, which is the basis of $\Lambda _{W}$ dual to $\Delta $.

Let $W$ be the Weyl group of $R$.  This is a finite group that acts on
$E$ via a reflection representation, and preserves the inner product
and the lattices $\Lambda _{R}$ and $\Lambda _{W}$.  There is a
decomposition of $E$ into a finite union of rational polyhedral cones,
and $W$ acts by permuting these cones.  Let $C^{+}$ be the closed
top-dimensional cone generated by $\Omega $.

Let $\Pi$ denote the set of isomorphism classes of complex irreducible
representations of $\ggg $.  It is known that elements of $\Pi $ are
in bijection with the set $\Lambda _{W}\cap C^{+}$, the dominant
weights.  Given $\lambda $ from this latter set, we denote the
corresponding representation by $\pi _{\lambda }$.  Then the
definition of the zeta function associated to $\ggg $ is
\begin{equation}\label{wzf}
\zeta _{\ggg } (s) := \sum _{\lambda \in \Lambda _{W}\cap C^{+}} (\dim
\pi _{\lambda })^{-s}.
\end{equation}

\subsection{}
Let $m>1$ be an integer.  The special value $\zeta _{\ggg } (2m)$
can be computed using a Dedekind sum as follows.  

Let $\rho$ be one-half the sum of the positive roots.  An application
of the Weyl character formula \cite[Corollary 24.6]{fulton.harris}
shows that for any dominant weight $\lambda $, we have
\begin{equation}\label{wdf}
\dim \pi _{\lambda } = \prod _{\alpha \in R^{+}} \frac{(\rho + \lambda
, \alpha )}{(\rho ,\alpha )}.
\end{equation}

It is known that any dominant weight $\lambda$ can be written as a
nonnegative integral linear combination of the fundamental weights.
Using this in \eqref{wdf}, a computation shows that \eqref{wzf}
becomes
\begin{equation}\label{wzf2}
\zeta _{\ggg } (2m) = M^{2m}\sum _{x\in (\Z ^{>0})^{\ell
}}\frac{1}{\prod _{i=1}^{r} \langle a_{i} , x\rangle ^{2m}}.
\end{equation}
Here $M$ is the integer $\prod _{\alpha \in R^{+}} (\rho ,\alpha )$,
$r$ is the number of positive roots, and the $a_{i} \in (\Z
^{>0})^{\ell}$ are the coefficients of the positive roots in terms of
$\Delta $.  The pairing $\langle \phantom{a},\phantom{a}\rangle$ is
the usual scalar product on $\R ^{\ell }$.

\subsection{}
We obtain a Dedekind sum by extending the sum \eqref{wzf2} to the
whole lattice.  Let $Z ^{\ell }\subset \R ^{r}$ be the span of the
first $\ell$ basis vectors, and let $\sigma = \sigma (\ggg ) $ be an
$r\times r$ integral matrix such that $\langle \sigma _{i},x\rangle =
\langle a_{i}, x\rangle $ for $x\in Z^{\ell }$, and such that $\det
\sigma = 1$.  Let $e = (2k, \dots ,2k)\in \R ^{r}$.

\begin{proposition}
(cf. \cite[p. 507]{zagier.zeta}) 
$\zeta _{\ggg } (2k) = \frac{M^{2k}}{\#W}S (Z^{\ell },\sigma (\ggg ),e,0)$.
\end{proposition}

Hence these special values can be computed in polynomial time using
our techniques.  We conclude with two examples:
$\SL_{3}$ and $\SL_{4}$.  We recommend verification of these formulas
to the interested reader for a pleasant combinatorial exercise.

\begin{proposition}\label{sln}
Let $\zeta (s)$ be the Riemann zeta function.  Then
\begin{align}
&\frac{6}{2^{2m}}\zeta _{\SL_{3}} (2m) = 8 \sum _{\substack{0\leq i\leq 2m\\i\equiv 0\mod{2}}}{\binom{4m-i-1}{2m-1}} \zeta (i)\zeta (6m-i)\label{sl3form}.\\
&\frac{24}{12^{2m}}\zeta _{\SL_{4}} (2m) = 16 \sum _{0\leq i\leq 2m}{\binom{4m-i-1}{2m-1}} (A+B+C+D),\quad \hbox{where}\\
A=\sum& _{\substack{0\leq j\leq 2m\\0\leq t\leq 4m+i-j\\j,t\equiv 0 \mod{2}}}\binom{2m+i-j-1}{i-1}\binom{6m+i-j-t-1}{2m-1}\zeta (j)\zeta (t)\zeta (12m-j-t),\\
B=\sum& _{\substack{0\leq j\leq 2m\\0\leq u\leq 2m\\j,u\equiv 0 \mod{2}}}\binom{2m+i-j-1}{i-1}\binom{6m+i-j-u-1}{4m+i-j-1}\zeta (j)\zeta (u)\zeta (12m-j-u),\\
C=\sum& _{\substack{0\leq k\leq i\\0\leq v\leq 4m+i-k\\k,v\equiv 0 \mod{2}}}\binom{2m+i-k-1}{i-k}\binom{6m+i-k-v-1}{2m-1}\zeta (k)\zeta (v)\zeta (12m-k-v),\\
D=\sum& _{\substack{0\leq k\leq i\\0\leq w\leq 2m\\k,w\equiv 0 \mod{2}}}\binom{2m+i-k-1}{i-k}\binom{6m+i-k-w-1}{4m+i-k-1}\zeta (k)\zeta (w)\zeta (12m-k-w).
\end{align}
\end{proposition}

\begin{remark}
The formula \eqref{sl3form} was independently discovered by Zagier,
S. Garoufalidis, and L. Weinstein \cite[p. 506]{zagier.zeta}.
\end{remark}

\begin{example}
Here are some special values of $\zeta _{\SL_{3}}$ and $\zeta _{\SL_{4}}$.

\begin{center}
\begin{tabular}[ht]{c||l}
$2m$&$(6m+1)! \cdot 6\cdot \zeta_{\SL _{3}} (2m)/ ( 2^{2m}\cdot (2\pi) ^{6m})$\\
\hline
\hline
2&$1/(2\cdot 3)$\\
4&$19/ (2\cdot 3\cdot 5)$\\
6&$1031/ (3\cdot 7)$\\
8&$(11\cdot 43\cdot 751)/ (2\cdot 7)$\\
10&$(5\cdot 13 \cdot 27739097)/(3\cdot 11)$\\
12&$(17\cdot 29835840687589)/ (3\cdot 5\cdot 7\cdot 13)$\\
14&$(2\cdot 17\cdot 19\cdot 89\cdot 127\cdot 6353243297)/ 7$\\
16&$(19\cdot 23\cdot 31\cdot 221137132669842886663)/ (2\cdot 5^{2}\cdot 13\cdot 17)$
\end{tabular}
\end{center}
\bigskip
\begin{center}
\begin{tabular}[ht]{c||l}
$2m$&$(12m+1)!\cdot (6m+1)\cdot  (4m+1)\cdot  24\cdot \zeta _{\SL _{4}} (2m)/ (12^{2m}\cdot (2\pi) ^{12m})$\\
\hline
\hline
$2$&$23/ 2$\\
$4$&$(3\cdot 7\cdot 14081)/ 2$\\
$6$&$(757409\cdot 23283173)/ (5\cdot 7)$\\
$8$&$(3\cdot 11\cdot 1021\cdot 5529809\cdot 754075957)/ 2$\\
$10$&$(13\cdot 116763209\cdot 1872391681\cdot 3187203549787)/ (5\cdot 11)$\\
$12$&$(17\cdot 1798397149\cdot 5509496891\cdot 6127205846988571484743)/ (3\cdot7\cdot 13)$
\end{tabular}
\end{center}
\end{example}

%
%

\bibliographystyle{amsplain}
\bibliography{dedekind}

\end{document}